\documentclass[11pt]{article}

\title{\sc\Large Orbit equivalence and actions of $\F_n$}
\author{Asger T\"ornquist}
\date{October 9, 2005}

\usepackage{amsmath}
\usepackage{amstext}
\usepackage{amssymb}
\usepackage{amsthm}
\usepackage{enumerate}

\DeclareMathOperator{\lh}{lh} \DeclareMathOperator{\inn}{Inn}
 \DeclareMathOperator{\Inn}{Inn}

\newcommand{\ftwo}{\mathbb F_2}
\newcommand{\mpt}{\mathcal M_\infty}

\newcommand{\psiab}{\psi_{\alpha,\beta}}
\newcommand{\psiba}{\psi_{\beta,\alpha}}

\DeclareMathOperator{\tfa}{TFA}

\def\R{{\mathbb R}}
\def\C{{\mathbb C}}

\def\N{{\mathbb N}}
\def\Z{{\mathbb Z}}

\def\F{{\mathbb F}}

\def\T{{\mathbb T}}

\newcommand{\sformat}[1]{\setcounter{paragraph}{0}\addtocounter{section}{1}\begin{center}{\sc\normalsize
 \Roman{section}. #1}\end{center}}

\begin{document}
\pagestyle{myheadings} \markboth{{\small {\it Orbit Equivalence
and $\F_n$}}}{{\small{Orbit Equivalence and actions of $\F_n$, A.
T\"ornquist}}}

\maketitle

\begin{abstract}
In this paper we show that there are ``$E_0$ many'' orbit
inequivalent free actions of the free groups $\mathbb F_n$, $2\leq
n\leq\infty$ by measure preserving transformations on a standard
Borel probability space. In particular, there are uncountably many
such actions.
\end{abstract}

\section*{\sformat{Introduction}}

Let $G_1,G_2$ be countable groups, acting by measure preserving
transformations on standard Borel probability spaces $(X_1,\mu_1)$
and $(X_2,\mu_2)$ respectively, giving rise to orbit equivalence
relations $E_{G_1}$ and $E_{G_2}$. We say that the actions of
$G_1$ and $G_2$ are \emph{orbit equivalent} if there is a measure
preserving bijection $\varphi:X_1\to X_2$ such that

$$
x E_{G_1} y \iff \varphi(x) E_{G_2} \varphi(y)
$$
almost everywhere.

The fundamental theorem in the study of the notion of orbit
equivalence is the theorem of H. Dye (\cite{dye},\cite{dye2}),
which states that two ergodic measure preserving actions of $\Z$
are orbit equivalent. Ornstein and Weiss (\cite{ornsteinweiss},
\cite{cofewe}) showed that this theorem extends to all countable
amenable groups.

The work of Connes, Weiss and Schmidt
(\cite{schmidt},\cite{connesweiss}), and more recently Hjorth
(\cite{hjorth}), shows that this characterizes amenability: A
countable group is amenable if and only if it has, up to orbit
equivalence, only one ergodic action by measure preserving
transformations on a standard Borel probability space.

\bigskip

The study of orbit equivalence is naturally related to the study
of Borel equivalence relations, and the notion of {\it Borel
reducibility}. (The reader may find a more thorough discussion of
Borel reducibility in \cite{jkl}.) For equivalence relations $E$
and $F$ on Polish spaces $X$ and $Y$, we say that $E$ is {\it
Borel reducible} to $F$, written $E\leq_B F$, if there is a Borel
function $f:X\to Y$ such that
$$
x E y\iff f(x) F f(y).
$$
In other words, $E$ is Borel reducible to $F$ if we can classify
points in $X$ up to $E$ equivalence by a Borel assignment of
invariants, that are $F$ equivalence classes.

An equivalence relation $E$ on a Polish space $X$ is said to be
\emph{smooth} or {\it concretely classifiable} if $E\leq_B
=_{2^\N}$, where $=_{2^\N}$ is the equality relation on the Cantor
space $2^\N$. In other words, a smooth equivalence relation admits
a Borel assignment of real numbers as complete invariants,
classifying elements of $X$ up to $E$ equivalence.

There are equivalence relations which are {\it not} smooth. The
cardinal example of such an equivalence relation is $E_0$, defined
on $2^\N$ by
$$
fE_0g \iff (\exists N)(\forall n\geq N) f(n)=g(n).
$$
It is not hard to see that $E_0\nleq_B =_{2^{\N}}$, and
$=_{2^{\N}}\leq_B E_0$, i.e. that $=_{2^{\N}}<_B E_0$. Hence we
cannot in a Borel way classify points up to $E_0$ equivalence
using real numbers as invariants. This can be understood as saying
that, in the sense of $\leq_B$, there are many more $E_0$ classes
than there are real numbers.

\bigskip

\noindent If $E_0\leq_B F$ for an equivalence relation $F$, we
will say that $F$ has (at least) ``$E_0$ many'' equivalence
classes. In this paper we show:

\paragraph{Theorem 1.} There are (at least) $E_0$ many orbit inequivalent
almost everywhere free actions of $\F_n$, $2\leq n\leq \infty$, by
measure preserving transformations on a standard Borel probability
space.

\bigskip

This may be seen as a strengthening of the following result of
Gaboriau and Popa:

\paragraph{Theorem (Gaboriau-Popa, \cite{gapo}.)} There are continuum many orbit
inequivalent a.e. free actions of $\F_n$, $2\leq n\leq\infty$, by
measure preserving transformations on a standard Borel probability
space.

\bigskip
It is worth noting explicitly that Theorem 1 is stronger than the
result of Gaboriau and Popa since it rules out the possibility of
finding a reasonable real-valued complete invariant for orbit
equivalence. In the light of Gaboriau's work on the notion of
``cost'' \cite{gab}, one might have hoped otherwise.

\paragraph{Outline and organization.} The proof of Theorem 1 relies mostly
on elementary methods, and does not involve the operator algebra
techniques used by Gaboriau and Popa in their result in
\cite{gapo}. That said, both results rely on similar ideas with
origins in Popa's work on rigidity phenomena and the notion of
relative property (T) in \cite{popa}.

We first consider in section II a particular (a.e.) free action of
$\F_2$ on a standard Borel probability space $(X,\mu)$ by measure
preserving transformations. This action has the special property,
that there is a countable group $G\subseteq L_0(X,\T)\subseteq
L^\infty(X)$, invariant under the action of $\F_2$, such that the
induced semi-direct product $G\rtimes\F_2$ has the relative
property (T).

Our strategy is to first obtain $E_0$ many actions of $\F_3$.
Denote by $\mpt(X)$ the (Polish) group of all measure preserving
transformations on $X$. In section III we prove a general lemma
which has a consequence that there is a dense $G_\delta$ set of
transformations that extends a given a.e. free m.p. action of
$\F_2$ to an a.e. free m.p. action of $\F_3$.

The main argument is presented in section V. We consider the
special action of $\F_2$ mentioned above, and introduce an
equivalence relation $\mathbf R$ on $\mpt(X)$ by letting $S\mathbf
R S'$ whenever the equivalence relation induced by the
transformation $S$ and (the action of) $\F_2$, and the equivalence
relation induced by $S'$ and $\F_2$, are orbit equivalent.
Similarly, we let $S \mathbf F S'$ whenever $S$ and $\F_2$ induce
the same equivalence relation as $S'$ and $\F_2$ (a.e.)

Using the relative Kazhdan property, we show that $\mathbf
R/\mathbf F$ has countable classes. It then follows easily that
$\mathbf R$ is meagre in $\mpt(X)\times\mpt(X)$. We then show that
a Theorem of Becker and Kechris applies to give us that $E_0\leq_B
\mathbf R$, which shows Theorem 1 in the case $\F_3$.

The case $\F_n$, $n>3$ is obtained similarly. At the end of
section V, we see that the case $\F_2$ follows from that of $\F_3$
by an expansion argument.

\bigskip

We also get the following surprising corollary:

\paragraph{Corollary.} Equality a.e. of equivalence relations
induced by a.e. free measure preserving actions of $\F_n$, $n\geq
2$, is not smooth.

\bigskip

The author does not know if either orbit equivalence or equality
in this context is {\it strictly} above $E_0$, however it seems
natural suspect that this is so.

\bigskip

The author is grateful to Sorin Popa for many helpful discussions
regarding this manuscript. The author also wishes to thank the
referee for his many useful comments and suggestions that have
improved this paper.

\section*{\sformat{An action of $\F_2$}}

In this section we will describe a particular ergodic action of
$\F_2$ which has some very interesting properties. The idea behind
this is due to Sorin Popa (\cite{popa}), and is fundamental to the
entire argument of this paper. First we recall the notion of
\emph{relative property} (T) and related concepts.

\paragraph{Definition (\cite{behava},\cite{hava}).} If $(\pi,\mathcal H)$ is unitary representation of
 a countable group $G$, $Q\subseteq G$ is a
subset, and $\varepsilon>0$, we say that a vector $v\in \mathcal
H$ is $(Q,\varepsilon)$-invariant if
$$
\sup_{g\in Q}
\|\pi(g)v-v\|_\mathcal{H}<\varepsilon\|v\|_\mathcal{H}.
$$
The semidirect product of countable groups $H\rtimes G$ has the
\emph{relative property} (T)\footnote{In the litterature one
sometimes refers to this situation by saying that {\it $(G,H)$ is
a pair with property} (T), cf. \cite{behava} 1.4.3.}, if there is
a finite set $Q\subseteq H\rtimes G$ and $\varepsilon>0$ such that
whenever $(\pi,\mathcal H)$ is a unitary representation of
$H\rtimes G$ with $(Q,\varepsilon)$-invariant vectors, then there
is a non-zero $H$-invariant vector.

\bigskip
If $H\rtimes G$ has the relative property (T), it can be seen (as
in \cite{behava} proposition 1.1.8) that given $\delta>0$ we may
find $\varepsilon>0$ such that if $v$ is a
$(Q,\varepsilon)$-invariant vector for $(\pi,\mathcal H)$, then
there is an $H$-invariant vector $v'$ such that
$\|v-v'\|_\mathcal{H}\leq \delta$. In other words, we can ensure
that almost invariant vectors are $\delta$ close to $H$-invariant
vectors.

\paragraph{Examples.}

The semidirect product $\Z^2\rtimes SL_2(\Z)$, corresponding to
the natural action of $SL_2(\Z)$ on $\Z^2$, has the relative
property (T). (\cite{kazhdan},\cite{burger}, \cite{shalom}).

The matrices
$$
A=\left(\begin{array}{cc} 1 & 2\\ 0 & 1\end{array}\right),
B=\left(\begin{array}{cc} 1 & 0\\ 2 & 1\end{array}\right)
$$
generate a copy of $\F_2$ inside $SL_2(\Z)$. It can be shown that
this subgroup has finite index, so that $\F_2$ is a lattice in
$SL_2(\Z)$. It follows (see \cite{behava} Theorem 1.5.1) that the
corresponding semidirect product $\Z^2\rtimes\F_2$ also has the
relative property $(T)$.

\paragraph{The action.}

We now describe the particular action of $\F_2$ with which we will
be working throughout this paper. We denote by $\T$ the 1-torus
$$
\T=\{z\in\C:|z|=1\}.
$$
Note that the map $x\mapsto e^{2\pi ix}$ identifies $\R/\Z$ and
with $\T$.

We consider the group $\Z^2$ and in particular it's dual, $\hat
\Z^2$. Every character of $\Z^2$ has the form
$$
\chi_{\mathbf{a}}\binom{n}{m}=e^{2\pi i(na_1+ma_2)}
$$
where $\mathbf{a}=\binom{a_1}{a_2}\in\R^2$. Hence we may identify
$\hat \Z^2$ with the 2-torus $\T^2$.

The action of $SL_2(\Z)$ on $\Z^2$ induces an action on
$\hat\Z^2$, as defined by
$$
\sigma\cdot\chi\binom{n}{m}=\chi(\sigma^{-1}\binom{n}{m}).
$$
We then have the formula
$$
\chi_{\mathbf{a}}=\chi_{(\sigma^{-1})^T\mathbf{a}},
$$
which shows that the corresponding action of $SL_2(\Z)$ on $\T^2$
is measure preserving. We have previously noted that $\F_2$ sits
inside $SL_2(\Z)$ as a lattice. The restriction of the action of
$SL_2(\Z)$ on $\T^2$ to this subgroup is the the action of $\F_2$
we wish to consider in detail.

The action of $SL_2(\Z)$ on $\T^2$ induces an action on
$L_0(\T^2,\T)\subseteq L^\infty(\T^2)$, where $L_0(\T^2,\T)$ is
the group of measurable functions $\T^2\to\T$ with pointwise
multiplication, through
$$
\sigma\cdot f(\chi)=f(\sigma^{-1}\cdot\chi).
$$
Define $\phi_{\binom{n}{m}}\in L_0(\T^2,\T)$ by
$$
\phi_{\binom{n}{m}}(\chi)= \chi\binom{n}{m}.
$$
Then the map $\Z^2\to L^\infty(\T^2)$ defined by
$$
\binom{n}{m}\mapsto \phi_{\binom{n}{m}}
$$
is an endomorphism of $\Z^2$ into $L_0(\T^2,\T)$, and we denote by
$G$ the group
$$
G=\{\phi_{\binom{n}{m}}:\binom{n}{m}\in \Z^2\}.
$$
It follows easily from the definitions that
$$
\sigma\cdot\phi_{\binom{n}{m}}=\phi_{\sigma\binom{n}{m}}, \
(\sigma\in SL_2(\Z)),
$$
which shows that $G$ is invariant under the action of $SL_2(\Z)$,
and the semidirect product $G\rtimes SL_2(\Z)$ is isomorphic to
$\Z^2\rtimes SL_2(\Z)$. It follows that $G\rtimes SL_2(\Z)$ has
the relative property (T). Note also that $G$ separates points in
$\T^2$.

\bigskip

Finally we observe:

\paragraph{Claim.} The action of the subgroup
$\F_2<SL_2(\Z)$ on $\T^2$ is ergodic.

\begin{proof} We must show that for any $f\in L_0^2(\T^2)=1^\perp$,
we have $\sigma\cdot f=f$ for all $\sigma\in\F_2$ iff $f=0$. It is
a standard fact of Fourier analysis on locally compact abelian
groups that $\{\phi_{\binom{n}{m}}:\binom{n}{m}\in\Z^2\}$ forms an
orthonormal basis for $L^2(\T^2)$ (see \cite{reiter}, p. 143 ff.).
Hence
$$
\{\phi_{\binom{n}{m}}:\binom{n}{m}\neq 0\}
$$
is an orthonormal basis for $L_0^2(\T^2)$.

Suppose first that $f\neq 0$ and
$$
f=a_1\phi_{\binom{n_1}{m_1}}+\cdots+a_l\phi_{\binom{n_l}{m_l}},
$$
$\binom{n_i}{m_i}\neq (0,0)$ and $a_1,\ldots, a_l\in\C$. Since
$$ A^k=\left(\begin{array}{cc} 1 & 2k\\ 0 &
1\end{array}\right), B^k=\left(\begin{array}{cc} 1 & 0\\ 2k &
1\end{array}\right)
$$
we may find $k>0$ such that for any $m_i\neq 0$ we have $(A^k
\binom{n_i}{m_i})_1\neq n_j$ for all $j\leq l$. Then find $k'>0$
such that for all $m_i=0$ we have $(B^{k'}\binom{n_i}{m_i})_2\neq
m_j$ for all $j\leq l$. Then for $\sigma=B^{k'}A^{k}$ we have
$\sigma\binom{n_i}{m_i}\neq \binom{n_j}{m_j}$ for all $j\leq l$,
and
$$
\sigma\cdot
f=a_1\phi_{\sigma\binom{n_1}{m_1}}+\cdots+a_l\phi_{\sigma\binom
{n_l}{m_l}},
$$
so that $(\sigma\cdot f,f)=0$. Now for general $0\neq f\in
L_0^2(\T^2)$ it follows that $|(\sigma\cdot f,f)|$ can be made
arbitrarily small for appropriate $\sigma\in\F_2$ and hence there
is $\sigma\in \F_2$ such that $\sigma\cdot f\neq f$.
\end{proof}

Let us summarize what we have shown:

\paragraph{Proposition 1.} There is an ergodic action of $\F_2$ on a
standard non-atomic probability space $(X,\mu)$, such that

\begin{enumerate}[(i)]

\item There is a countable group $G\subseteq L_0(X,\T)\subseteq
L^\infty(X)$, invariant under the induced action of $\F_2$ on
$L^\infty(X)$,

\item The group $G$ separates points in $X$, and

\item The semidirect product $G\rtimes\F_2$ has the relative
property (T).
\end{enumerate}

\bigskip

This is all we will need to know for the argument below.

\section*{\sformat{The Category Lemma}}

Let $(X,\mu)$ be a standard non-atomic probability space, and
denote by $\mpt(X)$ the group of all measure preserving
transformations on $X$. This group has two important group
topologies (see \cite{halmos} pp. 61 and 69): The first is the
\emph{weak topology} which is defined by the neighborhood basis
$$
N(T; E_1,\ldots, E_k, \varepsilon)=\{S\in\mpt(X) : (\forall i\leq
k) \mu(T(E_i)\triangle S(E_i)))<\varepsilon\}
$$
for $T\in\mpt(X)$, where $\varepsilon>0$ and $E_1,\ldots,
E_k\subseteq X$ are measurable subsets of $X$. With the weak
topology $\mpt(X)$ is a Polish group.

The other topology is the \emph{uniform topology}, which is
induced by the metric
$$
d_U(T,S)=\mu(\{x:T(x)\neq S(x)\}).
$$
The uniform topology is stronger than the weak topology and the
metric $d_U$ is complete. However, the uniform topology is not
separable.

The uniform topology will be useful in later sections, but for the
considerations of this section we only need the weak topology.

\bigskip

In this section we prove the following category theoretic fact:

\paragraph{The Category Lemma.} Let $G=\{T_n\in\mpt(X):n\in \N\}$ be a
countable group of measure preserving transformations, and suppose
$G$ acts freely almost everywhere on $(X,\mu)$. Then
$$
\{S\in\mpt(X):G * \langle S\rangle \text{ acts a.e. freely on
$X$}\}
$$
is a dense $G_\delta$ subset of $\mpt(X)$ in the weak topology.

\bigskip

Here $G *\langle S\rangle$ denotes the free product of the group
$G$ and the group generated by $S$, denoted $\langle S\rangle$,
which may be thought of formally as the set of finite sequences of
alternatingly elements from $G$ and elements of $\langle
S\rangle$, with the obvious concatenate-and-reduce operation as
composition.

\bigskip

Before the proof, let us note some useful facts:

\paragraph{Observation.}
If $(A_i)_{i\leq k}$ is a sequence of $k\in \N$ measurable subsets
of $(X,\mu)$, then there is a measure preserving involution
$P\in\mpt(X)$ such that $P(x)\neq x$ a.e. and $P(A_i)=A_i$ for all
$i\leq k$. The proof is an easy induction on $k$: For $k=1$, it is
simply the fact that there always is an involution with almost no
fixed points. Now suppose suppose the assertion holds for $k\geq
1$, and consider $(A_i)_{i\leq k+1}$. Then there is an involution
$T_0:A_{k+1}\to A_{k+1}$ with almost no fixed points such that
$T_0(A_i\cap A_{k+1})=A_i\cap A_{k+1}$, and similarly
$T_1:X\setminus A_{k+1}\to X\setminus A_{k+1}$. Then $T=T_0\cup
T_1$ is the desired transformation.

\bigskip

We also note the following easy technical lemma:

\paragraph{Lemma 2.} Let $T_1,\ldots, T_n:X\to X$ be measurable functions such that
$$
\mu(\{x : (\forall i,j) i\neq j\implies T_i(x)\neq T_j(x)\})>K>0.
$$
Then there are finitely many disjoint non-null Borel sets
$E_1,\ldots, E_m$ such that
$$
\mu(\bigcup_{l\leq m} E_l)>K,
$$
and $T_i(E_l)\cap T_j(E_l)=\emptyset$ whenever $i\neq j$, $l\leq
m$.

\begin{proof}
Assume that $X$ is equipped with a compatible Polish topology. By
Lusin's Theorem (cf. \cite{kechris}, 17.12), let
$$
F\subseteq \{x :
(\forall i,j) i\neq j\implies T_i(x)\neq T_j(x)\}
$$
be a closed set of measure $\mu(F)>K$, such that all of
$T_1,\ldots, T_n$ are continuous on $F$. Then for each $x\in F$,
there is a basic open set $O_x$ such that $T_i(O_x\cap F)\cap
T_j(O_x\cap F)=\emptyset$. Since the collection of sets $O_x$ is
countable, and $\mu(\bigcup_{x\in F} O_x\cap F)=\mu(F)$, there are
finitely many $O_{x_1},\ldots, O_{x_k}$, such that
$$
\mu(\bigcup_{l\leq k} O_{x_l}\cap F)>K.
$$
After possibly breaking the collection $(O_{x_l})_{l\leq k}$ into
disjoint pieces, we obtain a collection of disjoint Borel sets
$B_1,\ldots, B_m\subseteq X$ with the same properties. Now let
$E_l=B_l\cap F$.
\end{proof}

\begin{proof}[Proof of the Category Lemma.]
Let $\langle a\rangle$ be an infinite cyclic group with a single
generator $a$, and consider the free product $G *\langle a\rangle$
of $G=\{T_n:n\in \N\}$ and $\langle a\rangle$. The free product
may be thought of as consisting of words in the alphabet $\mathcal
A=\{T_n:n\in\N\}\cup\{a,a^{-1}\}$, reduced according to the rules
of the respective groups.

Given such a word $w$, the \emph{evaluation map}
$e_w:\mpt(X)\to\mpt(X)$ associated to $w$ is the map that
associates to a transformations $S\in\mpt(X)$ the transformation
obtained by replacing $a$ with $S$ in the word $w$. Note that
since $\mpt(X)$ is a topological group, the evaluations map $e_w$
is continuous.

We will show that for a non-trivial reduced word $w$ in the
alphabet $\mathcal A$,
$$
\{S\in \mpt(X) : e_w(S)(x)\neq x \text{ a.e.}\}
$$
is a dense $G_\delta$ set in $\mathcal M_\infty(X)$.

The proof goes by induction on the length of the word $w$. Assume
that the above holds for all non-trivial reduced words $\eta$ with
$\lh(\eta)<n$, and let $\lh(w)=n$. Let $\varepsilon>0$ be given,
and consider the set
\begin{equation*}
G_\varepsilon =\{S\in\mpt(X) : e_w(S)(x)\neq x \text{ on a set of
measure }
>1-\varepsilon\}.
\end{equation*}
It suffices to show this set is ({\it i}) open and ({\it ii})
dense.

\bigskip

({\it i}) {\it The set $G_\varepsilon$ is open}.

\bigskip

\noindent This will follow from:

\paragraph{Claim.} If $P\in \mpt(X)$ and
$$
\mu(\{x: P(x)\neq x\})>K>0,
$$
then there is a neighborhood $N\subseteq\mpt(X)$ of $P$ such that
$$
\mu(\{x: S(x)\neq x\})>K
$$
for all $S\in N$.

\begin{proof}
Let $\delta>0$ be such that
$$
\mu(\{x: P(x)\neq x\})>K+\delta.
$$
By lemma 2 applied to the identity transformation $I$ and $P$,
there are disjoint Borel sets $E_1,\ldots, E_m$, such that
$P(E_l)\cap E_l=\emptyset$ and $\mu(\bigcup E_l)>K+\delta$. Now
consider the neighborhood
$$
N=N(P; E_1,\ldots, E_m, \frac{\delta} m) =\{S: (\forall l\leq m)
\mu(P(E_l)\triangle S(E_l))<\frac \delta m\}.
$$
Then for any $S\in N$ we have $\mu(S(E_l)\cap E_l)<\frac\delta m$,
so that
$$
\mu(\{x: S(x)\neq x\})\geq \sum_{l=1}^m \mu(E_l)-\frac \delta m >
K,
$$
and $N$ is the neighborhood we needed to find.
\end{proof}

From the claim it follows easily that $G_\varepsilon$ is open,
since
$$
G_{\varepsilon} = e_w^{-1}(\{P: P(x)\neq x \text{ on a set of
measure } >1-\varepsilon\}),
$$
where $e_w$ is the evaluation map.

\bigskip

({\it ii}) {\it The set $G_\varepsilon$ is dense}.

\bigskip

Fix a sequence $\eta_0,\ldots, \eta_n=w$ of reduced words in the
alphabet $\mathcal A$, such that $\lh(\eta_i)=i$, and for all
$i<n$ there is a unique $\tau\in\mathcal A$, such that $
\eta_{i+1}=\tau\eta_i$.

Let $S\in\mpt(X)$ and let $N$ be a neighborhood of $S$. We want to
show that there is $S'\in N\cap G_\varepsilon$. By our inductive
assumption, we can assume that
$$
S\in \bigcap_{\lh(\eta)<n}\{S': e_\eta(S')(x)\neq x \text{
a.e.}\}.
$$

Moreover, we can assume the letter $a$ (or $a^{-1}$) occurs at
some point in the reduced word $w$, since otherwise there is
nothing to show.

Clearly $e_{\eta_i}(S)(x)\neq e_{\eta_j}(S)(x)$ a.e.
whenever $i\neq j$, ($i,j<n$). So by lemma 2, we can find disjoint
non-0 measurable sets $E_1,\ldots, E_M\subseteq X$ such that
$$
\mu(\bigcup_{l\leq M} E_l)>1-\varepsilon
$$
and $e_{\eta_i}(S)(E_l)\cap e_{\eta_j}(S)(E_l)=\emptyset$ whenever
$i\neq j$.

Let $(A_i)_{i\in\mathcal F_0}$ be a finite family of measurable
sets and let $\delta>0$ be such that
$$
N_0 = N(S; (A_i)_{i\in\mathcal F_0}, \delta)\subseteq N,
$$
and such that the collection $(A_i)_{i\in\mathcal F_0}$ contains all of the sets
$$
e_{\eta_i}(S)(E_l),\ (l\leq M, i<n).
$$

Define
$$
f(P)=\{x: e_w(P)(x)=x\}
$$
for $P\in \mpt(X)$. Then either $\mu(f(S))<\varepsilon$ (in which
case there is nothing to show), or there is some $E_l$ with
$\mu(f(S)\cap E_l)>0$. We may assume that $\mu(f(S)\cap E_1)>0$.
Let $i_0<n$ be largest possible such that
$\eta_{i_0+1}=a\eta_{i_0}$ or $\eta_{i_0+1}=a^{-1}\eta_{i_0}$.
Without loss of generality, assume $\eta_{i_0+1}=a\eta_{i_0}$ and
let $B=e_{\eta_{i_0}}(S)(f(S)\cap E_1)$. As observed before the
proof, we may find an involution $P\in\mpt(X)$ such that $P(x)\neq
x$ for almost all $x\in B$, and $P(x)=x$ for all $x\notin B$, such
that $P(A_i)=A_i$ for all $i\in\mathcal F_0$.

Define
$$
S_1=SP.
$$
Then $S_1(A_i)=S(A_i)$ for all $i\in\mathcal F_0$, in particular
$S_1\in N_0$. Moreover, for almost all $x\in E_1\cap f(S)$, we
have $e_w(S_1)(x)\neq e_w(S)(x)=x$, and for $x\in E_1\setminus
f(S)$ we have $e_w(S_1)(x)=e_w(S)(x)\neq x$. Thus $e_w(S_1)(x)\neq
x$ for almost all $x\in E_1$.

By Lemma 2, we may then find disjoint measurable sets $F_1,\ldots
F_p\subseteq E_1$, such that
$$
\mu(\bigcup_{q\leq p} F_q\cup\bigcup_{1<l\leq M}
E_p)>1-\varepsilon
$$
and $e_w(S_1)(F_q)\cap F_q=\emptyset$, $q\leq p$.

Let $(A_i)_{i\in \mathcal F_1}$ be the extension of the family
$(A_i)_{i\in \mathcal F_0}$ obtained by adding all the sets
$$
e_{\eta_i}(S_1)(F_q), \ (i<n, q\leq p).
$$

If $\mu(f(S_1))<\varepsilon$, then we're done. Otherwise we may
find $l>1$ for which $\mu (f(S_1)\cap E_l)>0$, indeed we may
assume $\mu (f(S_1)\cap E_2)>0$. Now we can apply the same
argument as above, with $(A_i)_{i\in \mathcal F_1}$ in place of
$(A_i)_{i\in\mathcal F_0}$, to get $S_2\in N_0$ with $f(S_2)\cap
E_2=\emptyset$, and a finite collection $(A_i)_{i\in\mathcal F_2}$
extending $(A_i)_{i\in\mathcal F_1}$. However, as this
construction guarantees that $S_1(A_i)=S_2(A_i)$ for all $i\in
\mathcal F_1$, we retain that $e_w(S_2)(F_q)\cap F_q=\emptyset$,
for $q\leq p$. Hence by repeating the above argument finitely many
times, we eventually obtain a transformation $S'=S_l\in N_0$, for
some $l\leq M$, such that
$$
\mu(\{x : e_w(S')(x)\neq x\})
>1-\varepsilon,
$$
and since $S'\in N_0\subseteq N$, this completes the proof.
\end{proof}

\paragraph{Remark.} The author is thankful to Sorin Popa for making
the following remarks regarding the Category Lemma in a recent
conversation.

The proof of the category lemma does not use the fact that
$G=\{T_n\in\mpt(X):n\in\N\}$ is a group. Rather, the proof shows
that given a sequence of transformations
$\{T_n\in\mpt(X):n\in\N\}$ with almost no fixed points, we can
find a dense $G_\delta$ set of $T\in\mpt(X)$ that are
``independent'' of $\{T_n:n\in\N\}$, in the sense that any
composition of alternatingly elements $T_n$, $n\in\N$, and $T$ or
$T^{-1}$, does not have any fixed points almost everywhere. Popa
has pointed out that this observation gives us the following
useful corollary:

\paragraph{Corollary.} Let $(X,\mu)$ be a standard Borel probability
space, and suppose $H_1$ and $H_2$ are countable groups of
transformations in $\mpt(X)$, acting freely a.e. on $X$. Then
there is a dense $G_\delta$ set of transformations $T\in\mpt(X)$
such that $H_1*TH_2T^{-1}$ acts freely a.e. on $X$.
\begin{proof}
By the remark, we can find a dense $G_\delta$ set of
transformations $T\in\mpt(X)$ that are independent of the set of
transformations $H_1\cup H_2$, in the above sense. But then
$H_1*TH_2T^{-1}$ acts freely a.e. on $X$.
\end{proof}

\section*{\sformat{The group $\mpt(X)$}}

Before we proceed to prove Theorem 1, we note in this section some
standard facts regarding the group $\mpt(X)$ and some important
subgroups.

If $G$ is a countable group acting by m.p. transformations on the
standard Borel probability space $(X,\mu)$, giving rise to the
equivalence relation $E_G$, the \emph{full group}, or \emph{inner
group}, of $E_G$ is the group
$$
\inn(E_G)=\{S\in\mpt(X) : S(x)E_G x \text{ a.e.}\}.
$$
It is easy to see that in this case $\inn(E_G)$ is a Polish group
when given the \emph{uniform} topology. The full group is also
denoted $[E_G]$. The set of partial measure preserving functions,
i.e. \emph{morphisms}, whose graph is contained in $E_G$, is
denoted $[[E_G]]$.

\paragraph{Proposition 2.} Let $G$ be a countable group acting by m.p. transformations on the
standard Borel probability space $(X,\mu)$. Consider $\mpt(X)$
with the weak topology. Then

\begin{enumerate}[(i)]
\item $\inn(E_G)$ is a meagre subgroup of $\mpt(X)$, and

\item $\inn(E_G)$ is dense if and only if $E_G$ is ergodic.
\end{enumerate}
\begin{proof}
(i) Since the uniform topology is stronger than the weak topology,
the identity embeds $\Inn(E_G)$ continuously into $\mpt(X)$, and
so $\Inn(E_G)$ is an analytic (in fact, Borel) subgroup of
$\mpt(X)$. In particular, $\Inn(E_G)$ has the Baire property, and
it follows by Pettis Theorem (cf. \cite{kechris}) that either it
is meagre, or it contains a neighborhood of the identity. However,
the latter cannot be the case, since then $\Inn(E_G)$ would not be
separable in the uniform topology.

\bigskip

(ii) The ``only if'' direction is obvious. For the ``if''
direction, let $A,A'\subseteq X$ be measurable sets with
$\mu(A)=\mu(A')>0$. We first show, that there is a morphism
$\varphi:A\to A'$, $\varphi\in[[E_G]]$. Since
$$
\mu(\bigcup_{g\in G} g\cdot A)=1,
$$
there is some $g\in G$ such that $\mu(g\cdot A\cap A')>0$. Let
$B_0=g^{-1}\cdot A'\cap A$ and define $\varphi(x)=g\cdot x$ for
$x\in B_0$. If $\mu(A\setminus B_0)=0$, we're done. Otherwise, we
can repeat the argument with $A\setminus B_0$ and $A'\setminus
g\cdot B_0$. In this way, we eventually exhaust $A$ (in perhaps
transfintely many steps), and have defined the desired morphism
$\varphi$.

Let $T\in\mpt(X)$, and let $A_1,\ldots, A_k$ be measurable subsets
of $X$. We claim that there is $S\in\inn(E_G)$ such that
$\mu(S(A_i)\triangle T(A_i))=0$ for all $i\leq k$. After possibly
breaking the sets $A_1,\ldots, A_k$ into smaller pieces, we can
assume they are disjoint. Then by the above we can find morphisms
$\varphi_0,\ldots,\varphi_k\in [[E_G]]$ such that
$\varphi_i:A_i\to T(A_i)$ for $1\leq i\leq k$ and
$$
\varphi_0: X\setminus \bigcup A_i\to T(X\setminus \bigcup A_i).
$$
Define $S=\bigcup_{i\geq 0} \varphi_i$.

Note now, that $S\in N(T;A_1,\ldots,A_k,\varepsilon)$ for all
$\varepsilon>0$. Since $T$ and $A_1,\ldots,A_k$ were arbitrary,
this shows that $\inn(E_G)$ is dense in $\mpt(X)$.
\end{proof}

Let $(O_n)$ be a countable basis for a compatible Polish topology
on the standard probability space $(X,\mu)$. Assume the basis
$(O_n)$ is closed under finite unions. Then a complete metric for
the weak topology on $\mpt(X)$ is given by
$$
d_w(S,T)=\sum_{m} 2^{-m}\big[\mu(S(O_m)\triangle
T(O_m))+\mu(S^{-1}(O_m)\triangle T^{-1}(O_m))\big],
$$
(see \cite{kechris}.) We note the following fact about the
relation between convergence in $d_w$ and pointwise convergence:

\paragraph{Proposition 3.} Let $S\in\mpt(X)$ and suppose $(S_n)$ is a sequence of measure preserving
transformations such that $d_w(S_n,S)<2^{-n}$. Then $S_n(x)\to
S(x)$ a.e.
\begin{proof}
Let $\varepsilon>0$ and $\rho>0$. Let $F\subseteq X$ be a closed
set such that $S$ is continuous on $F$ and $\mu(F)>1-\frac \rho
2$. We want to show that
$$
\mu(\{x: (\exists N)(\forall n\geq N) \
d(S_n(x),S(x))<\varepsilon\})>1-\rho.
$$
For this, first find finitely many basic open sets
$O_{m_1},\ldots,O_{m_k}$ such that
$$
x,y\in O_{m_i}\cap F\implies d(S(x),S(y))<\varepsilon
$$
and $\mu(\bigcup_{i\leq k} O_{m_i}\cap F)>1-\frac \rho 2$.

Since $d_w(S_n,S)<2^{-n}$, we have for each $i$ that
$\mu(S_n(O_{m_i})\triangle S(O_{m_i}))<C_{m_i}2^{-n}$, where
$C_{m_i}>0$ is a constant which depends only on $m_i$. Then for
$N>0$ such that $C_{m_i}2^{-N}<\frac \rho {2k}$,
$$
\mu(S(O_{m_i})\cap \bigcap_{n> N}
S_n(O_{m_i}))>\mu(S(O_{m_i}))-\frac \rho {2k}.
$$
Since if $x,y\in S(O_{m_i}\cap F)$ we have $d(x,y)<\varepsilon$,
it now follows that
$$
\mu(\{x: (\exists N)(\forall n\geq N) \
d(S_n(x),S(x))<\varepsilon\})>1-\rho.
$$
\end{proof}

\section*{\sformat{The main argument}}
In this section we prove Theorem 1. We focus on proving the
theorem for $\F_3$, the argument for $\F_n$, $n> 3$, being
similar. The case $\F_2$ will eventually follow from that of
$\F_3$.

Consider an action of $\F_2=\langle a,b\rangle$ on a standard
probability space $(X,\mu)$, as described in Proposition 1. Denote
by $G\subseteq L_0(X,\T)\subseteq L^{\infty}(X)$ the associated
$\F_2$-invariant multiplicative subgroup, and let $T_a(x)=a\cdot
x$ and $T_b(x)=b\cdot x$ be the m.p. transformations corresponding
to the generators $a$ and $b$. Let $Q\subseteq G\rtimes\F_2$ be a
finite Kazhdan set with Kazhdan constant $\varepsilon>0$,
witnessing the relative property (T) of $G\rtimes\F_2$.

For $S\in\mpt(X)$ denote by
$$
E_S=E_{\langle T_a,T_b,S\rangle}
$$
the equivalence relation generated by the transformations
$T_a,T_b$ and $S$. We define two equivalence relations $\mathbf R$
and $\mathbf F$ on $\mpt(X)$ by
$$
S\mathbf RS'\iff E_S\text{ is orbit equivalent to } E_{S'}
$$
and
$$
S\mathbf FS'\iff E_S=E_{S'} \text{ a.e.}
$$
Denote by $A\subseteq \mpt(X)$ the set of transformations $S$ such
that $\langle T_a,T_b, S\rangle$ induces an a.e. free action of
$\F_3$ on $X$. It follows from the Category Lemma that this set is
a dense $G_\delta$ set. Then Theorem 1 in the case of $\F_3$ can
be phrased as

\paragraph{Theorem 1$^\prime$.} $E_0\leq_B\mathbf R|A$.

\bigskip

An outline of the proof is as follows: We will first show that
$\mathbf F$ has meagre classes. Then we will use the relative
property (T) to make an argument modeled on \cite{popa} and
\cite{hjorth}, to show that $\mathbf R$ has countable index over
$\mathbf F$, and deduce that it is a meagre subset of
$\mpt(X)\times\mpt(X)$. It will then be easy to apply a theorem of
Becker and Kechris (\cite{beke}) to obtain that $E_0\leq_B \mathbf
R|A$.

The proof is presented as a sequence of lemmata. We start by
computing the complexity of $\mathbf R$ and $\mathbf F$:

\paragraph{Lemma 3.} The equivalence relations $\mathbf R$ and
$\mathbf F$ are analytic.

\begin{proof} Give the space $X$ a compatible Polish topology,
and let furthermore the metric $d_w$ on $\mpt(X)$ be as in
Proposition 3. Let $(S_n)_{n\in\mathbb N}$ be a sequence of Borel
measure preserving transformations which is dense in $\mpt(X)$.
Define a relation $\tilde E(\phi, x,y)\subseteq \mathbb N^{\mathbb
N}\times X\times X$ by
$$
\tilde E(\phi, x,y)\iff S_{\phi(n)}(x)\to y.
$$
Then $\tilde E$ is Borel. Define $\Phi:\mpt(X)\to \mathbb
N^{\mathbb N}$ by letting $\Phi(S)(n)$ be the least $m\in \mathbb
N$ such that
$$
d_w(S,S_m)<2^{-n}.
$$
Then clearly $d_w(S_{\Phi(S)(n)},S)<2^{-n}$ for all $n$, and
$\Phi$ is a Borel map. Hence the set $E\subseteq \mpt(X)\times
X\times X$ defined by
$$
E(S,x,y)\iff \tilde E(\Phi(S),x,y)
$$
is Borel. By Proposition 3,
$$
(\forall^{\mu}x,y) [E(S,x,y) \iff S(x)=y]
$$
and
\begin{align*}
S\mathbf F S' &\iff (\forall^\mu x,y) [xE_Sy\iff xE_{S'} y]\\
&\iff (\forall^\mu x,y) [(\exists \tau) e_\tau(S)(x)=y\iff
(\exists\tau') e_{\tau'}(S')(x)=y]\\
&\iff (\forall^{\mu} x,y)[(\exists \tau) E(e_\tau(S),x,y)\iff
(\exists\tau') E(e_{\tau'}(S'),x,y)],
\end{align*}
where $\tau,\tau'$ are words in $\{a,b,c\}$, and $e_\tau$ denotes
the evaluation map, as in section III (i.e. $T_a$ is substituted
for $a$, $T_b$ is substituted for $b$, and $S$ is substituted for
$c$). Since the measure quantifiers preserve analyticity (see
\cite{kechris}, p. 233), we conclude that $\mathbf F$ is analytic
(in fact, $\mathbf F$ is easily seen to be Borel.)

Finally,
$$
S\mathbf R S'\iff  (\exists T\in\mpt(X)) (\forall^{\mu}x,y) [x E_S
y\iff T(x) E_{S'} T(y)]\\
$$
So that $\mathbf R$ is analytic, since
\begin{align*}
(\forall^\mu x,y) \big[T(x) E_{S'} T(y)
\iff & (\forall z,z') E(T,x,z)\wedge E(T,y,z')\implies zE_{S'} z'\\
\iff & (\exists z,z') E(T,x,z)\wedge E(T,y,z')\wedge zE_{S'}
z'\big]
\end{align*}
\end{proof}

\paragraph{Corollary 1.} The equivalence relation $\mathbf F$ is
meagre in $\mpt(X)\times\mpt(X)$, and has meagre classes.

\begin{proof}
Since $\mathbf F$ is analytic, it has the Baire property. Hence by
\cite{kechris} 8.41, it is enough to show that each $\mathbf
F$-class is meagre. But $[S]_{\mathbf F}\subseteq\inn(E_S)$, which
is meagre by Proposition 2 (i).
\end{proof}

\bigskip

Our next step is to show

\paragraph{Main Lemma.} The equivalence relation $\mathbf R/\mathbf F$ is countable. That is, each $\mathbf R$-class contains at most countably many $\mathbf F$-classes.

\bigskip

Before the proof, we note the following:

\paragraph{Observation.} It is an easy observation, that if $(Y,d)$ is a Polish space, and
$(y_\alpha)_{\alpha<\omega_1}$, is a sequence, then for every
$\delta>0$ there is an unbounded set $B\subseteq\omega_1$ such
that whenever $\alpha,\beta\in B$, then
$d(y_\alpha,y_\beta)<\delta$.

Similarly, if $G$ is a Polish group and
$(g_\alpha)_{\alpha<\omega_1}$ is a sequence in $G$, then for any
neighborhood $N\subseteq G$ of the identity in $G$, there is an
unbounded set $B\subseteq \omega_1$ such that whenever
$\alpha,\beta\in B$, then $g_\alpha g_\beta^{-1}\in N$. To see
this, associate to any $g\in G$ a basic open neighborhood $N_g$ of
$g$, such that for $(h_1,h_2)\in N_g\times N_g$, we have
$h_1h_2^{-1}\in N$, using the continuity of the group operations.
Since $(N_g)_{g\in G}$ is countable there must be some $g_0$ such
that $g_\alpha\in N_{g_0}$ for an unbounded set of $\alpha$.

\bigskip

We also note:

\paragraph{Lemma 4.} For $g\in L^\infty(X)$ and $\delta>0$, there is
a neighborhood $N$ of the identity $I\in\mpt(X)$, such that
$$
\psi\in N\implies \| g-g\circ\psi\|_{L^2(X)}<\delta.
$$

\begin{proof}
This is trivial if we note that the weak topology on $\mpt(X)$ is
precisely the subspace topology inherited from the unitary group
on $L^2(X)$, under the identification $\psi\mapsto U_\psi$ where
$U_\psi(g)=g\circ\psi^{-1}$, $g\in L^2(X)$.
\end{proof}

Finally, recall that if $G$ is a countable group acting by m.p.
transformations on $(X,\mu)$, a Borel measure $M$ on $E_G$ is
defined by
$$
M(A)=\int |A_x|d\mu(x)
$$
for $A\subseteq E_G$ (see \cite{kemi}, p. 34).

\begin{proof}[Proof of Main Lemma.] Let $(S_\alpha)_{\alpha<\omega_1}$ be a sequence of
m.p. transformations such that $S_\alpha \mathbf R S_\beta$ for
all $\alpha,\beta\in \omega_1$. We want to show that
$E_{S_\alpha}=E_{S_\beta}$ for some $\alpha\neq\beta$.

For each $\alpha<\omega_1$, let $\psi_\alpha:X\to X$ be a m.p.
transformation witnessing that $E_{S_0}$ is orbit equivalent to
$E_{S_\alpha}$. We will write $\psiab$ for the transformation
$\psi_\beta\circ\psi_\alpha^{-1}$. A unitary representation
$\pi_{\alpha,\beta}$ of $G\rtimes \ftwo$ is defined on
$L^2(E_{S_\alpha})$ by
$$
((g,\sigma)\cdot
f)(x,y)=g(x)\overline{g(\psiab(y))}f(\sigma^{-1}\cdot x,
\psiba\sigma^{-1}\psiab(y))
$$
for each $\alpha<\omega_1$. Indeed, we have
\begin{align*}
&((g_1,\sigma_1)\cdot((g_2,\sigma_2)\cdot f))(x,y)\\
=&g_1(x)\overline{g_1(\psiab(y))}((g_2,\sigma_2)\cdot
f)(\sigma_1^{-1}\cdot x,
\psiba\sigma_1^{-1}\psiab(y))\\
=&g_1(x)\overline{g_1(\psiab(y))}g_2(\sigma_1^{-1}\cdot
x)\overline{g_2(\sigma_1^{-1}\psiab(y))}\\
&f(\sigma_2^{-1}\sigma_1^{-1}\cdot x,
\psiba\sigma_2^{-1}\sigma_1^{-1}\psiab(y))\\
=&((g_1(\sigma_1\cdot g_2),\sigma_1\sigma_2)\cdot f))(x,y),
\end{align*}
which shows that an action of $G\rtimes \F_2$ is defined.

\paragraph{Claim.} There is $\alpha\neq \beta$ such that
$\psiab(x)=x$ on a non-null set of $x$.

\begin{proof}
Recall that $G\rtimes \ftwo$ has the relative property (T) with
Kazhdan pair $(Q,\varepsilon)$. We assume that $Q$ has the form
$Q_1\times Q_2$, for finite sets $Q_1\subseteq G$ and
$Q_2\subseteq \ftwo$, and that $\varepsilon$ is chosen so that if
$(\pi,\mathcal H)$ is a unitary representation and $v$ is
$(Q,\varepsilon)$-invariant, then there is a $G$-invariant $v'$
such that $\|v-v'\|_{\mathcal H}\leq \frac 1 2$.

In order to prove the claim, we will show, that there is
$\alpha\neq\beta$ such that the function $\mathbf 1_\Delta\in
L^2(E_{S_\alpha})$,
$$
\mathbf 1_\Delta(x,y)=\left\{\begin{array}{ll}

1 & \text{ if $x=y$}\\

0 & \text{ otherwise}
\end{array}\right.
$$
is $(Q,\varepsilon)$ invariant for the representation
$\pi_{\alpha,\beta}$. If we can show this, then there is a
$G$-invariant $\xi\in L^2(E_{S_\alpha})$ such that $\|\xi-\mathbf
1_{\Delta}\|_{L^2(E_{S_\alpha})}\leq\frac 1 2$. From this it
follows that $\xi(x,x)\neq 0$ on a non-null set and hence
$$
g(x)\overline{g(\psiab(x))}=1\text{ for all $g\in G$}
$$
on a non-null set. Since $G$ separates points, it follows that
$\psiab(x)=x$ on a non-null set, as we wanted.

Since $Q_1$ is finite, we can use Lemma 4 to find a neighborhood
$N\subseteq\mpt(X)$ of $I\in\mpt(X)$ such that for each $\psi\in
N$ and $g\in Q_1$ we have
$$
\|g-g\circ\psi\|_{L^2(X)}^2<\frac {\varepsilon^2} {2}.
$$
Using the observation preceeding the proof, we find an unbounded
set $B_0\subseteq \omega_1$ such that $\psiab\in N$ for all
$\alpha,\beta\in B_0$.

We now consider the transformations
$\psi_\alpha^{-1}\sigma^{-1}\psi_\alpha\in\inn(E_{S_0})$,
$\sigma\in\ftwo$. Applying the first part of the observation, we
can for a given $\sigma\in Q_2$ find an unbounded set
$B_1\subseteq B_0$ such that
\begin{equation}
d_U(\psi_\alpha^{-1}\sigma^{-1}\psi_\alpha,
\psi_\beta^{-1}\sigma^{-1}\psi_\beta)<\frac {\varepsilon^2}
{4}\tag{*}
\end{equation}
for all $\alpha,\beta\in B_1$, where $d_U$ is the usual complete
metric for the uniform topology on $\inn(E_{S_0})$. Iterating this
until the finite set $Q_2$ is exhausted, we get an unbounded set
$B\subseteq B_0$ such that (*) holds for all $\sigma\in Q_2$ and
all $\alpha,\beta\in B$.

By (*) it holds for $\alpha,\beta\in B$ that the set
$$
C_{\sigma}=\{x\in X: \psiba\sigma^{-1}\psiab(x)=\sigma^{-1}\cdot
x\}
$$
has $\mu(C_{\sigma})>1-\frac {\varepsilon^2} {4}$ for $\sigma\in
Q_2$.

Consider then the unitary representation $\pi_{\alpha,\beta}$ for
some fixed $\alpha,\beta\in B$. For $(g,\sigma)\in Q=Q_1\times
Q_2$, we have
\begin{align*} &\|\mathbf 1_\Delta\ -\pi_{\alpha,\beta}(g,\sigma)
\mathbf 1_\Delta\|^2_{L^2(E_{S_\alpha})}\\
=& \int\sum_{y\in[x]_{E_{S_\alpha}}}| \mathbf 1_\Delta(x,y)-
g(x)\overline{g(\psiab(y))}\mathbf 1_\Delta(\sigma^{-1}\cdot x,
\psiba\sigma^{-1}\psiab(y))|^2 d\mu(x).
\end{align*}
Since for almost all $x$
\begin{align*}
\sum_{y\in[x]_{E_{S_\alpha}}} & | \mathbf 1_\Delta(x,y)-
g(x)\overline{g(\psiab(y))}\mathbf 1_\Delta(\sigma^{-1}\cdot x,
\psiba\sigma^{-1}\psiab(y))|^2\\
&\leq  1+\|g\|_\infty^2\|\overline g\|_\infty^2=2,
\end{align*}
we get
\begin{align*} &\|\mathbf 1_\Delta\ -\pi_{\alpha,\beta}(g,\sigma)
\mathbf 1_\Delta\|^2_{L^2(E_{S_\alpha})}\\
&\leq \int_{C_{\sigma}}\sum_{y\in[x]_{E_{S_\alpha}}}| \mathbf
1_\Delta(x,y)- g(x)\overline{g(\psiab(y))}\mathbf
1_\Delta(\sigma^{-1}\cdot x,
\psiba\sigma^{-1}\psiab(y))|^2d\mu(x)\\
&+ \frac {\varepsilon^2} 2\\
&\leq \int_{C_{\sigma}} | \mathbf 1_\Delta(x,x)-
g(x)\overline{g(\psiab(x))}\mathbf 1_\Delta(\sigma^{-1}\cdot x,
\psiba\sigma^{-1}\psiab(x))|^2d\mu(x) + \frac{\varepsilon^2} 2\\
&\leq \| 1-g(\overline{g\circ
\psiab})\|^2+\frac {\varepsilon^2} 2\\
&\leq\|\overline{g\circ \psiab}\|^2 \| g\circ
\psiab-g\|^2+\frac {\varepsilon^2} 2\\
&\leq \varepsilon^2.
\end{align*}
Hence  $\mathbf 1_\Delta$ is $(Q,\varepsilon)$-invariant, as
claimed.
\end{proof}

Let $\alpha\neq \beta$ as in the claim. Since $E_{\ftwo}$ is
$\mu$-ergodic we may assume, after possibly discarding a set of
measure zero, that in each $E_{\ftwo}$ class there is $x$ such
that
$$
\psiab(x)=x.
$$

Consider an $E_{S_\alpha}$ class $C=[x]_{E_{S_\alpha}}$. We claim
that $C\subseteq [x]_{E_{S_\beta}}$. For this, write
$$
C=\bigcup_{x_i} [x_i]_{E_{\ftwo}}
$$
where each $x_i$ is such that
$$
\psiab(x_i)=x_i.
$$
Then $x_i E_{S_\beta} x_j$ for all $i,j$, and so
$$
C=\bigcup_{x_i} [x_i]_{E_{\ftwo}}\subseteq [x]_{E_{S_\beta}}.
$$
The opposite inclusion follows by a similar argument, and we
conclude that $E_{S_\alpha}=E_{S_\beta}$.
\end{proof}

\paragraph{Corollary 2.} The relation $\mathbf R$ is meagre in
$\mpt(X)\times\mpt(X)$.

\begin{proof}
By Corollary 1 and the Main Lemma, the $\mathbf R$-classes are
meagre, so we can conclude that $\mathbf R$ is meagre as in the
proof of Corollary 1.
\end{proof}

\bigskip

Recall, that $A$ is the set of transformations $S$, such that
$\langle T_a,T_b, S\rangle$ act freely a.e. on $X$. We now prove:

\paragraph{Theorem 1$^\prime$.} $E_0\leq_B\mathbf R|A$.

\begin{proof} We will use the following theorem:

\paragraph{Theorem (Becker-Kechris, \cite{beke} proof of 3.4.5, also \cite{hjorth} p. 32.)}
Suppose $E$ is an equivalence relation on the Polish space $X$,
which is meagre as a subset of $X\times X$.  Suppose further, that
there is a group $G$ acting by homeomorphisms on $X$, such that
$E_G\subseteq E$, and that there is a dense $G$-orbit. Then
$E_0\leq_B E$.

\bigskip

Let $G=\inn(E_{\mathbb F_2})$ act on $\mpt(X)$ by conjugation. For
$S\in\mpt(X)$ and $T\in G$, it is clear that
$$
E_{TST^{-1}}\subseteq E_{S}.
$$
But for the same reason
$$
E_S=E_{T^{-1}TST^{-1}T}\subseteq E_{TST^{-1}},
$$
so that $E_S=E_{TST^-1}$. Hence $E_G\subseteq \mathbf F\subseteq
\mathbf R$, and $G$ acts by homeomorphisms.

For an aperiodic transformation $S\in\mpt(X)$, the conjugacy class
of $S$ in $\mpt(X)$ is dense (see \cite{halmos}, p. 77). Since by
Proposition 2, $G=\inn E_{\mathbb F_2}$ is dense in $\mpt(X)$, it
follows that $[S]_G$ is dense in $\mpt$ for any aperiodic $S$. In
particular, $E_G$ has a dense orbit.

Since the set $A$ is easily seen to be invariant under the action
of $G$ it now follows from Becker-Kechris' theorem that $E_0\leq_B
\mathbf R|A$.
\end{proof}

\paragraph{Remark.} It is worth noting that the above proof gives us as a corollary that $E_0\leq_B \mathbf F$.
In particular, \emph{equality}
a.e. of equivalence relations induced by
actions by m.p. transformations is \emph{not} concretely
classifiable (i.e. smooth).

\paragraph{The case $\F_2$.} As remarked earlier, a virtually
identical argument to the above can be made for the case $\F_n$,
$3<n\leq\infty$: We simply add the appropriate number of
generating, independent transformations.

The case $\F_2$ follows from the case $\F_3$ by an expansion
argument. Consider again the family $E_S=E_{\langle
T_a,T_b,S\rangle}$, $S\in A$, of equivalence relations on
$(X,\mu)$, as above. Let $Y=X\times\{0,1\}$ with the product
measure (with equal weight to $0$ and $1$.) Let
$\tau(x,i)=(x,1-i)$, which is measure preserving, and define for
each $S\in\mpt(X)$ the transformation $\tilde S\in\mpt(Y)$ by
$\tilde S(x,0)=(S(x),0)$, and $\tilde S(x,1)=(x,1)$. For each
$S\in A$, let $\tilde E_S$ be the equivalence relation on $Y$
generated by $\tilde T_a, \tilde T_b, \tilde S$ and $\tau$. If
$E_S$ is orbit equivalent to $E_{S'}$, then clearly $\tilde E_S$
is orbit equivalent to $\tilde E_{S'}$. Conversely, suppose
$\tilde\varphi\in\mpt(Y)$ witnesses that $\tilde E_{S}$ is orbit
equivalent to $\tilde E_{S'}$. Define for $i,j\in\{0,1\}$
$$
C_{ij}=\{x\in X: \tilde\varphi(x,i)=(y,j) \text{ for some $y\in
X$}\}.
$$
Then $\mu(C_{01})=\mu(C_{10})$. By Proposition 2 (ii), there is a
morphism $\psi\in[[E_S]]$ such that $\psi(C_{01})=C_{10}$. Now for
$x\in C_{00}$, define $\varphi(x)=y$ where
$\tilde\varphi(x,0)=(y,0)$, and for $x\in C_{01}$, define
$\varphi(x)=y$ where $\tilde\varphi(\psi(x),1)=(y,0)$. Then
$\varphi$ witness that $E_S$ is orbit equivalent to $E_{S'}$.

Finally, the equivalence relation $\tilde E_S$, $S\in A$, is
generated by the transformations
$$
T_0(x,i)=\left\{\begin{array}{ll}
\tau(x,i) & \text{ if $i=0$,}\\
\tilde T_a \tau(x,i) & \text{ if $i=1$,}
\end{array}\right.
$$
and
$$
T_1(x,i)=\left\{\begin{array}{ll}
\tau\tilde S(x,i) & \text{ if $i=0$,}\\
\tilde T_b \tau(x,i) & \text{ if $i=1$.}
\end{array}\right.
$$
which induce an a.e. free m.p. action of $\F_2$. This shows that
there are $E_0$ many orbit inequivalent actions of $\F_2$, and
this completes the proof of Theorem 1.

\paragraph{Remark.} The author does not know if orbit equivalence
for free m.p. actions of $\F_n,n\geq 2$, is in general strictly
more complicated than $E_0$. In the case of Kazhdan groups, it
turns out that there are at least $\tfa$ many orbit inequivalent
actions, where $\tfa$ denotes the isomorphism relation for
countable torsion free abelian groups, cf. \cite{atphd}. In
particular, it follows from a result of Hjorth \cite{hjorth3} that
orbit equivalence is analytic non-Borel in this case. Hence the
author finds it natural to suspect that in the case $\F_n, n\geq
2$, orbit equivalence is also far more complicated than $E_0$.

\bigskip

\begin{small}
{\sc\noindent Department of Mathematics\\
University of California at Los Angeles\\
Los Angeles, CA 90095-1555\\
{\it E-mail}: atornqui@math.ucla.edu}
\end{small}

\end{document}